\documentclass[titlepage,11pt]{article}
\usepackage{latexsym}
\usepackage{amsfonts}
\usepackage{amsmath}
\newcommand\blackslug{\hbox{\hskip 1pt \vrule width 4pt height 8pt depth 1.5pt
        \hskip 1pt}}
\newcommand\bbox{\hfill \quad \blackslug \bigbreak}

\def\DD{\hbox{-}}

\def\LL{,\ldots,}

\title{Thomassen's theorem on the two-linkage problem in acyclic digraphs: a shorter proof}
\author{
Paul Seymour\thanks{Supported by AFOSR grant
FA9550-22-1-0234, and NSF grant  DDMS-2154169.}\\
Princeton University, Princeton, NJ 08544}

\date{September 14, 2024; revised \today}

\newtheorem{thm}{}[section]

\newcommand{\Proof}{\noindent{\bf Proof.}\ \ }

\begin{document}
\maketitle
\begin{abstract}
Let $G$ be an acyclic digraph, and let $a,b,c,d\in V(G)$, where $a,b$ are sources, $c,d$ are sinks, and every other 
vertex has in-degree and out-degree at least two. In 1985, Thomassen showed that there do not exist disjoint directed paths 
from $a$ to $c$
and from $b$ to $d$, if and only if $G$ can be drawn in a closed disc with $a,b,c,d$ drawn in the boundary in order. 
We give a shorter proof.

\end{abstract}

\section{Introduction}
Digraphs in this paper are finite, and without loops or parallel edges. A digraph is {\em acyclic} if it has no directed cycle, and 
a {\em dipath} is a directed path.
The {\em $k$-linkage problem} is, given two $k$-tuples $(a_1\LL a_k)$ and $(b_1\LL b_k)$ of vertices of a digraph $G$, 
to decide whether there are $k$ pairwise vertex-disjoint dipaths $P_1\LL P_k$ in $G$, where $P_i$ is from $a_i$ to $b_i$ 
for each $i$.
The $k$-linkage problem is NP-complete~\cite{fortune} for general digraphs $G$, even if $k=2$, but solvable in polynomial time~\cite{fortune} if $G$ is acyclic.

The 2-linkage problem in acyclic digraphs is particularly nice, because it comes equipped with a theorem characterizing the digraphs 
in which the problem is infeasible. It is easy to reduce the characterization question to the case handled by the following beautiful
theorem of Thomassen~\cite{thomassen}:
\begin{thm}\label{thomassenthm}
Let $G$ be an acyclic digraph, and let $S=\{a,b\}$ be its set of sources, and $T=\{c,d\}$ its set of sinks. 
Suppose that $S\cap T=\emptyset$, and every vertex not in $S\cup T$ has in-degree and out-degree both at least two. Then exactly one of the following holds:
\begin{itemize}
\item there are vertex-disjoint dipaths $P,Q$ of $G$, where $P$ is from $a$ to $c$ and $Q$ is from $b$ to $d$;
\item $G$ can be drawn in a closed disc with $a,b,c,d$ drawn in the boundary in order.
\end{itemize}
\end{thm}
Thomassen's proof was about five pages long. Our objective is to give a shorter proof of a slightly stronger statement.

Let $G$ be a digraph, and let $S=\{s_1\LL s_k\}$ and $T=\{t_1\LL t_{\ell}\}$ be disjoint subsets of $V(G)$. Let $P$ be a dipath from 
$s_i\in S$ to $t_j\in T$, and let $P'$ be a dipath from $s_{i'}\in S$ to $t_{j'}\in T$. We say the pair
$(P,P')$ is a {\em cross
relative to $(s_1\LL s_k), (t_1\LL t_{\ell})$ } if $P,P'$ are vertex-disjoint 
and either $i<i'$ and $j>j'$, or $i>i'$ and $j<j'$. 
We will prove the following (\ref{thomassenthm} is the case when $k=\ell=2$):
\begin{thm}\label{mainthm}
Let $G$ be an acyclic digraph, and let $S=\{s_1\LL s_k\}$ be its set of sources, and $T=\{t_1\LL t_{\ell}\}$ its set of sinks.
Suppose that $S\cap T=\emptyset$, and every vertex not in $S\cup T$ has in-degree and out-degree both at least two. Then exactly one of the following holds;
\begin{itemize}
\item there is a cross in $G$ relative to $(s_1\LL s_k), (t_1\LL t_{\ell})$;
\item $G$ can be drawn in a closed disc with $s_1,s_2\LL s_k, t_{\ell}, t_{\ell-1}\LL t_1$ drawn in the boundary in order.
\end{itemize}
\end{thm}
\Proof
That not both statements hold is clear; so we assume that the first bullet is false, and will prove that the second is true
by induction on $|V(G)|+|E(G)|$. Let $C= V(G)\setminus (S\cup T)$. The following observation will be useful:
\\
\\
(1) {\em For all $u,v\in C\cup T$ (possibly equal), there are two dipaths from $\{u,v\}$ to $T$, vertex-disjoint if $u\ne v$, 
and vertex-disjoint except for $v$ if $u=v$. The same holds for all $u,v\in S\cup C$ and dipaths from $S$ to $\{u,v\}$.
}
\\
\\
Let $W=\{u\}\cap \{v\}$. If the first statement is false, then by Menger's theorem, there is a set $X\subseteq V(G)\setminus W$ 
with $|X|\le 1$ such that every dipath from $\{u,v\}$
to $T$ has a vertex in $X$. Not both $u,v\in X$ (since $X\cap W=\emptyset$), so we assume that $v\notin X$.
Choose a maximal dipath $P$ with first vertex $v$ and with no vertex in $X$, and let $p$ be
the last vertex of $P$. Thus $p\notin S$ since $v\notin S$ and $S$ is the set of sources; and $p\notin T$ from the choice of $X$. So
$p\in C$, and therefore has out-degree at least two; and so has an out-neighbour $q\notin X$. Since $G$ is acyclic, 
$q\notin V(P)$, and so we could add $q$ to $P$ and make a longer dipath, contrary to the maximality of $P$. This proves the first statement,
and the second follows similarly. This proves (1).
\\
\\
(2) {\em We may assume that there is no edge between $S,T$.}
\\
\\
Suppose that $s_i$ is adjacent to $t_j$, where $1\le i\le k$ and $1\le j\le \ell$. 
Let $A$ be the set of 
$v\in C$ such that there is a dipath from $\{s_1\LL s_{i-1}\}$ to $v$, and 
let $A'$ be the set such that there is a dipath from $\{s_{i+1}\LL s_k\}$ to $v$. Similarly, let $B$ be the 
set of $v\in C$ such that there is a dipath from $v$ to $\{t_1\LL t_{j-1}\}$, and let $B'$ be the set of $v\in C$ 
such that there is a dipath from $v$ to $\{t_{j+1}\LL t_{\ell}\}$. From (1), $A\cup A'=B\cup B'=C$.
But if $v\in A$, then $v\notin B'$, since otherwise the union of a dipath from $\{s_1\LL s_{i-1}\}$ to $v$ and one from 
$v$ to $\{t_{j+1}\LL t_{\ell}\}$ is a dipath that makes a cross with the dipath $s_i\DD t_j$. Hence $A\cap B'=\emptyset$, and similarly
$A'\cap B=\emptyset$; so $A=B$ and $A'=B'$. Define $V=A\cup \{s_1\LL s_{i-1}\} \cup \{t_1\LL t_{j-1}\}$, and $V'= 
A'\cup \{s_{i+1}\LL s_{k}\} \cup \{t_{j+1}\LL t_{\ell}\}$. 

Thus $V_1,V_2,\{s_i,t_j\}$ are pairwise disjoint and have union $V(G)$.
There is no edge from $V$ to $V'$,
since if $uv$ is such an edge then $u\notin T$ and $v\notin S$, and the union of a dipath from $\{s_1\LL s_{i-1}\}$ to $u$,
the edge $uv$, and a 
dipath from $v$ to $\{t_{j+1}\LL t_{\ell}\}$ is a dipath that makes a cross with the dipath $s_i\DD t_j$. Similarly there
is no edge from $V'$ to $V$. Let $H$ be obtained from $G$ by deleting the edge $s_it_j$.
The result follows from the inductive hypothesis, applied to $H\setminus V'$ with the sequences
$(s_1\LL s_i)$, $(t_j\LL t_{\ell})$ and applied to $H\setminus V$ with the sequences $(s_{i}\LL s_k), (t_1\LL t_{j})$.
This proves (2).

\bigskip

If $C=\emptyset$, then $E(G)=\emptyset$ by (2) and the result is true, so we assume that $C\ne \emptyset$. Hence we may choose $v\in C$ with no in-neighbour in $C$, since $G$ is acyclic. But $v$ has at least two in-neighbours in $G$,
and they all belong to $S$; let the in-neighbours of $v$ be $\{s_i:i\in I\}$, where $I\subseteq \{1\LL k\}$. Thus $|I|\ge 2$; let $h,j$
be the smallest and largest members of $I$. 
\\
\\
(3) {\em For $h<i<j$, $s_i$ has no out-neighbour except possibly $v$.}
\\
\\
Suppose that $s_iu$ is an edge where $u\ne v$. 
Hence $u\in C\cup T$, and by (1) there are two vertex-disjoint dipaths $P,Q$ from $\{u,v\}$ to $T$. 
One of $P,Q$, say $Q$, has first vertex $v$, and so $P$ has first vertex $u$. Let $Q_1$ be obtained from $Q$ by adding the edges $s_hv$, 
let $Q_2$ be obtained from $Q$ by adding $s_jv$, and let $P'$ be obtained from $P$ by adding $s_iu$. Then one of $(P',Q_1), (P',Q_2)$
is a cross relative to $(s_1\LL s_k), (t_1\LL t_{\ell})$ (which one depends on the order of the ends of $P,Q$ in $T$), a contradiction. 
This proves (3).

\bigskip

Let $G'$ be obtained from $G$ by deleting $s_{h+1}\LL s_{j-1}$ and the edges $s_hv, s_jv$, and let 
$S'=(S\setminus \{s_{h+1}\LL s_{j-1}\})\cup \{v\}$. Thus $S'\cap T=\emptyset$, and from (3), $S'$ is the set of sources of $G'$, and every vertex of
$G'$ not in $S'\cup T$ has out-degree and in-degree at least two. Suppose that in $G'$ there is a cross $(P,Q)$ 
relative to $(s_1\LL s_h,v,s_j\LL s_k), (t_1\LL t_{\ell})$. One of $P,Q$ has first vertex $v$ (since there is no cross in $G$
relative to $(s_1\LL s_k), (t_1\LL t_{\ell})$), say $P$. Since $h\ne j$, one of $s_h,s_j\notin V(Q)$, say $s_i$.
Let $P'$ be obtained by adding $s_iv$ to $P$. Then $(P',Q)$ is a cross in $G$
relative to $(s_1\LL s_k), (t_1\LL t_{\ell})$, a contradiction. So there is no such cross. From the inductive hypothesis,
the result holds for $G'$ and the sequences $(s_1\LL s_h,v,s_j\LL s_k), (t_1\LL t_{\ell})$; but then it holds for $G$, by (3). This proves \ref{mainthm}.~\bbox

\end{document}